\documentstyle[12pt]{article}

\newcommand{\sect}[1]{\setcounter{equation}{0}\section{#1}}

\def\be{\begin{equation}}
\def\ee{\end{equation}}
\def\bea{\begin{eqnarray}}
\def\eea{\end{eqnarray}}

\def\l{\lambda}
\def\m{\mu}
\def\n{\nu}

\def\jj{J_3}
\def\jp{J_+}
\def\jm{J_-}

\def\>#1{{\bf #1}}

\def\1{\'{\i}}
\def\R{{\rm I\kern-.2em R}}

\begin{document}

\begin{center} {\LARGE{\bf{Jordanian $U_{h,s}gl(2)$ and its }}}

 {\LARGE{\bf{coloured realization}}}
 \end{center}

\bigskip\bigskip\bigskip


\begin{center} Preeti Parashar\footnote{Present address: Departamento de
Fisica,
Universidad de Burgos, Burgos, Spain and Departamento de Fisica Teorica,
Universidad Autonoma de Madrid, Madrid, Spain.}
\end{center}

\begin{center} {\it  Department of Mathematics \\
Texas A  \& M
University,  College Station,
Texas 77843 , USA.}
\end{center}

\bigskip\bigskip\bigskip

\begin{abstract}
A two-parametric non-standard (Jordanian) deformation of the
Lie algebra $gl(2)$ is constructed, and then, exploited to obtain
a new, triangular $R$-matrix solution of the coloured Yang-Baxter
equation. The corresponding coloured quantum group is presented
explicitly.

\end{abstract}
\bigskip
\bigskip
\bigskip
\centerline{To appear in Lett. Math. Phys.}

\newpage


Non-standard or Jordanian  deformation of the group $GL(2)$ was
introduced in \cite{russ,zak}. Its two parameter generalization was given
in
\cite{aga} and the supersymmetric version $GL(1/1)$ was worked out
in \cite{dp}.
On the quantum algebra level , the non-standard deformation of the Lie
algebra $sl(2)$  was first proposed by Ohn \cite{ohn} .
The universal $R$-matrix for $U_hsl(2)$ was written down in
\cite{vla,iran,bur} and the irreducible representations were
studied in \cite{dob,chakra}.

A unique feature governing this deformation
is that the $R$-matrix which solves the quantum Yang-Baxter equation
(YBE) is triangular i.e. $R_{12}R_{21} =
I$.
This motivates us to ask a question : Is it possible to Yang-Baxterize
this solution such that  now the $R$-matrix  satisfies the more general
spectral parameter dependent YBE ?
It is quite obvious that the triangularity property poses the
main obstacle to carry out the usual Baxterization procedure. However,
a possible way out is , to make the $R$-matrix dependent on some
continuously varying colour parameters so that it obeys the
coloured YBE. The colour parameters play the same role as spectral
parameters  and the coloured YBE and spectral YBE become equivalent.

The coloured $R$-matrix solution was introduced by
Murakami \cite{mur} and various other approaches [12 - 18]
have been persued to construct them.
In an interesting work \cite{bh}, Burd\'ik and Hellinger  showed that
coloured solutions can be obtained from deformations of non-semisimple
Lie algebras like $U_qgl(2)$ which may be thought of as a splice product
of $U_qsl(2)$ and $Uu(1)$. The eigenvalues of the extra Casimir operator
coming from the invariant subalgebra $Uu(1)$ were interpreted as the
colour indices.
In the present letter, we adopt this approach to find a
{\it new , non-standard solution of the coloured } YBE corresponding
to $GL(2)$.
However, in order to carry out the analysis an immediate
problem arises as the non-standard deformation of $gl(2)$ is not known.
Hence, our purpose becomes two-fold. The first is, to give such a
quantization of
$gl(2)$,and the second is, to derive coloured solution from it. So let us
proceed with our first aim . We shall adhere to the convenient
basis of \cite{bur}.

\sect{The quantum algebra $U_{h,s}gl(2)$}

The classical Lie algebra $gl(2)$ is defined as
\be
[\jj,\jp]=2\jp,\qquad  [\jj,\jm]=-2\jm,\qquad [\jp,\jm]=\jj  ,
\qquad [Z,\cdot\,]=0,
\label{ga}
\ee
where $Z$ is the central generator of $u(1)$.

The non-standard classical $r$-matrix associated with $gl(2)$ can be
written as
a combination of non-standard $r$ for $sl(2)$ and a term containing the
two primitive generators $\jp$ and $Z$.
\be
r= h\jj\wedge\jp + s Z\wedge\jp .
\label{gb}
\ee
It is a solution of the classical Yang--Baxter equation  and generates the
Lie bialgebra with cocommutators given by
$\delta(X) = [1\otimes X+ X\otimes 1 , r ]$:
\bea
&&\delta(\jp)=0 , \qquad \delta(Z)=0 , \cr
&&\delta(\jj)=2h\, \jj\wedge \jp +2s\, Z\wedge \jp , \cr
&&\delta(\jm)=2h\, \jm\wedge \jp + s\, \jj\wedge Z.
\eea
where $h$ and $s$ are the two deformation parameters.

The Jordanian quantum algebra $U_{h,s}gl(2) = U_hsl(2) \oplus Uu(1)$ can
be defined as a Hopf algebra with the following structure:

Commutation relations:
\bea
&&[\jj,\jp ]= \frac{e^{2 h \jp} -1  } h,\qquad
[\jj,\jm]=-2 \jm +h \jj^2 +2 s Z \jj + {s^2\over h}Z^2,\cr
&&[\jp  ,\jm ]= \jj +{s\over h}Z(1 - e^{2h\jp}),\qquad [Z,\cdot\,]=0.
\label{gk}
\eea

Coproduct:
\bea
&&  \Delta (\jp) =1 \otimes \jp  + \jp \otimes 1,\cr
&& \Delta (\jj) =1 \otimes \jj  + \jj\otimes
e^{2h \jp } +{s\over h}Z \otimes (e^{2h \jp}- 1), \label{gj} \\
 && \Delta (\jm) = 1 \otimes \jm + \jm\otimes e^{2 h \jp } + s(
\jj + {s\over h}Z)\otimes Z\, e^{2 h \jp }, \cr
&&  \Delta (Z) =1 \otimes Z  + Z \otimes 1,\nonumber
\eea

Counit:
\be
\epsilon(X) =0,\qquad
 \mbox{for $X\in \{\jj,\jp,\jm,Z\}$},
 \ee

Antipode:
\bea
&&\gamma(\jp)=-\jp,\qquad  \gamma(\jj)=-\jj
e^{-2h\jp} + {s\over h}Z(1 - e^{-2h \jp}), \cr
&&\gamma(Z)=-Z,\qquad \gamma(\jm)=-\jm e^{-2h\jp} + s( \jj + {s\over h}Z )
Z e^{-2h\jp},
\eea

It can be checked that
$\{\jj ,\jp ,Z\}$ generates a Hopf
subalgebra of this.
Both the deformation parameters $h$ and $s$ play equal role in
shaping the above non-standard deformation,
unlike the two-parameter standard $(p,q)$-deformation \cite{swz}
where the second parameter is relegated to the coalgebra.
When $h = s$ ,this yields a one parameter
non-standard quantization of $gl(2)$.

\sect{Coloured $R$-matrix $R_{h,s}^{(\l,\m)}$}

The universal $R$-matrix for $U_{h,s}gl(2)$ is given by
\be
{\cal  R}=\exp\{- \jp\otimes (h \jj + s Z)\}\exp\{(h \jj + s Z)\otimes
\jp\} .
\ee
This element solves the quantum Yang--Baxter equation
\be
{\cal  R}_{12}{\cal  R}_{13}{\cal  R}_{23}=
{\cal  R}_{23}{\cal  R}_{13}{\cal  R}_{12}.
\label{aac}
\ee
and verifies the following quasi-triangular properties
$(X \in U_{h,s}gl(2))$:

\be
\sigma \circ \Delta( X) ={ \cal  R} \Delta( X)
{\cal  R}^{-1},  \qquad
\sigma(x\otimes y)=(y\otimes x),
\label{aaa}
\ee

\be
(\Delta  \otimes id){\cal  R}
={\cal  R}_{13}{\cal  R}_{23},\qquad
  (id \otimes \Delta){\cal  R} ={\cal  R}_{13}{\cal  R}_{12}.
\label{aab}
\ee

In addition, ${\cal R}^{-1}={\cal
R}_{21}$ (i.e. $\cal R$ is triangular ), and therefore $U_{h,s}gl(2)$
is endowed with a
triangular Hopf algebra structure.

Now let us denote the eigenvalue of $Z$ (Casimir like operator) by $\eta$
and the corresponding $n$-dimensional irreducible representation of the
quantum algebra (\ref{gk}) by $\pi_{\eta}^n$. Acting $\pi_{\eta}^n$
on the universal $\cal R$ yields a finite dimensional coloured matrix
representation . So
taking the following fundamental two-dimensional matrix representation of
$U_{h,s}gl(2)$
\bea
&& \pi_{\eta}(\jj)=\left(\begin{array}{ll}
 1 & 0  \\
 0 & -1 \end{array}\right),\quad
\pi_{\eta}(Z)=\left(\begin{array}{ll}
\eta & 0 \\
0 & \eta \end{array}\right),\cr
&& \pi_{\eta}(\jp)=\left(\begin{array}{cc}
 0 & 1  \\
0 & 0 \end{array}\right),\quad
\pi_{\eta}(\jm)=\left(\begin{array}{cc}
 (h+\eta s)^2/{2h} & 0  \\
 1 & (h-\eta s)^2/{2h} \end{array}\right),\label{ao}
\eea
and substituting in ($\eta$ acts as $\l$ in first factor of the tensor
product and $\m$ in the second )
\be
R^{(\l,\m)} = (\pi_{\l} \otimes \pi_{\m}) \cal R
\ee
we obtain the $4 \times 4$ matrix form of $\cal R$
\be
R^{(\lambda, \mu)}_{h,s}=
\left(\begin{array}{rrrr}
 1 & h+\lambda s & -(h+\mu s)& h^2-\lambda \mu s^2 -hs(\lambda-\mu)  \\
 0 & 1 & 0& h-\mu s   \\
 0 & 0 & 1 & -(h-\lambda s)  \\
 0 & 0 & 0& 1    \end{array}\right)
\label{ap}
 \ee
 We shall call this the coloured Jordanian $R$-matrix
, as it satisfies the coloured Yang- Baxter equation
\be
R^{(\l,\m)}_{12} R^{(\l,\n)}_{13} R^{(\m,\n)}_{23} =
R^{(\m,\n)}_{23} R^{(\l,\n)}_{13} R^{(\l,\m)}_{12}.
\label{aq}
\ee
$\l, \m, \n$ are referred to as the colour parameters, which vary
continuously. Setting $\l = \m = \eta$, and making the following
identification
\be
h+\eta s = z^{\prime}, ~~~~~ h- \eta s = z
\label{ar}
\ee
the coloured R-matrix (\ref{ap}) reduces to the two parameter $(z,
z^{\prime})$ $R$-matrix classified in \cite{hla}, and
in the limit  $\l = \m =0$ to the single parameter $(z)$ Jordanian
$R$-matrix \cite{russ}.

Analogous to the quantum group case,
one can also define the coloured Jordanian braid group representation
( BGR ) as
\be
{\hat R}^{(\l,\m)} = P R^{(\l,\m)}
\ee
where $R^{(\l,\m)}$ is taken as in (\ref{ap}) and $P$ is the usual
permutation
matrix. This coloured BGR turns out to be a solution of the
coloured braided YBE
\be
{\hat R}^{(\l,\m)}_{23} {\hat R}^{(\l,\n)}_{12} {\hat R}^{(\m,\n)}_{23} =
{\hat R}^{(\m,\n)}_{12} {\hat R}^{(\l,\n)}_{23} {\hat R}^{(\l,\m)}_{12}.
\label{at}
\ee
It is worth mentioning, that
${\hat R}^{(\l,\m)}$ does not satisfy the Hecke condition. Rather,
its characteristic equation is of fourth order given by
\be
{({\hat R}^{(\l,\m)} - 1 )}^3 ({\hat R}^{(\l,\m)} + 1 ) = 0, \l \ne \m
\ee
with only two distinct eigenvalues $+1$ and $-1$.

\sect{Coloured quantum group $GL_{h,s}^{(\l,\m)}(2)$}

Let us now investigate the quantum group structure associated with
the coloured Jordanian $R$-matrix.
The $RTT$ relations which normally define a quantum group
are modified suitably to
incorporate the coloured extension as \cite{frt}
\be
R^{(\l,\m)} T_{1\l} T_{2\m} = T_{2\m} T_{1\l} R^{(\l,\m)}
\label{baa}
\ee

To be compatible with the coloured $R$-matrix, now the
$T$-matrix  is also parametrized by the colour  $\l$ and
$\m$ (one at a time) .Since $\l, \m$ are continuous variables, this
implies the coloured quantum group has
an infinite number of generators.

As an illustration , let us consider the two-dimensional case.
Substituting
$R^{(\l,\m)}$ from (\ref{ap}) and $T_{\l}$ as ( and likewise $T_{\m}$)
\be
T_{\l}=\left(\begin{array}{cc}
 a_{\l} & b_{\l}  \\
 c_{\l} & d_{\l} \end{array}\right),\label{bab}
\ee
in the coloured RTT (\ref{baa}), we obtain the following independent
commutation relations for
the coloured Jordanian quantum group denoted by $GL_{h,s}^{(\l,\m)}(2)$:
\bea
&&[a_{\l},c_{\m}] = -(h-\m s)c_{\m}c_{\l}, \qquad
[a_{\l},d_{\m}] = (h+\l s)c_{\m}a_{\l} - (h-\m s)c_{\l}d_{\m}, \cr
&&[c_{\l},d_{\m}] = (h+\l s) c_{\m}c_{\l}, \qquad
[b_{\l},c_{\m}] = -(h+\m s)c_{\m}a_{\l} - (h-\m s)d_{\l}c_{\m}, \cr
&&[a_{\l},b_{\m}] = (h+\l s)a_{\m}a_{\l} -(h+\l s)a_{\l}d_{\m} +
(h+\m s)c_{\l}b_{\m} -f(\l,\m)c_{\l}d_{\m},\cr
&&[d_{\l},b_{\m}]=(h-\l s)d_{\l}d_{\m} -(h-\l s)a_{\m}d_{\l}
+(h-\m s)b_{\m}c_{\l} + f(\l,\m)a_{\m}c_{\l},\cr
&&[a_{\l},a_{\m}]=-(h+\l s)a_{\l}c_{\m} +(h+\m s)c_{\l}a_{\m} -
f(\l,\m)c_{\m}c_{\l},\cr
&&[b_{\l},b_{\m}]=-(h- \l s)a_{\m}b_{\l} +(h-\m s)b_{\m}a_{\l}
-(h+\l s)b_{\l}d_{\m} +(h+\m s)d_{\l}b_{\m}\cr
&&\qquad \qquad +f(\l,\m)(a_{\m}a_{\l} - d_{\l}d_{\m}),\cr
&&[c_{\l},c_{\m}] = 0,\cr
&&[d_{\l},d_{\m}]= -(h-\l s)c_{\m}d_{\l} + (h-\m s)d_{\m}c_{\l}
+ f(\l,\m)c_{\m}c_{\l}.
\label{bac}
\eea
Here, and in what follows

$f(\l,\m) =  h^2-\l\m s^2 -hs(\l-\m) ,~~~~
f(\m,\l) = h^2 - \l\m s^2 -hs(\m-\l)$.

All other relations can be generated from the set (\ref{bac}) by simply
interchanging $\l$ and $\m$. Notice the presence of the last four
relations like $[a_{\l},a_{\m}]$ etc. This is a peculier characteristic
of only the coloured quantum groups and does not exist for their
'uncoloured' counterparts since there, either $\l = \m$ or both
colours vanish.

\noindent
The coproduct and counit for the coloured algebra are given by
\cite{frt}
\be
\Delta (T_{\l}) = T_{\l} \otimes T_{\l}, \qquad
\epsilon (T_{\l}) = {\bf 1} .
\ee
Contrary to the algebra, the coalgebra is a function of only one colour
parameter.

\noindent
The quantum determinant corresponding to $T_{\l}$ is
\bea
&&{\cal D}_{\l} = a_{\l}d_{\l} - b_{\l}c_{\l} - (h+\l s)a_{\l}c_{\l}\cr
&& \qquad = a_{\l}d_{\l} - c_{\l}b_{\l} + (h-\l s)c_{\l}d_{\l}
\label{ada}
\eea
It can be checked that ${\cal D}_{\l}$ is still 'group like'
\be
{\cal D}_{\l}(T_{\l}{T_{\l}}^{\prime}) = {\cal D}_{\l}(T_{\l}).
{\cal D}_{\l}({T_{\l}}^{\prime}) ,\qquad if~~ [(T_{\l})_{ij}, ({
T_{\l}}^{\prime})_{kl}] = 0,
\ee
The  commutation relations of ${\cal D}_{\l}$ with the elements of
$T_{\l}$ are:
\bea
&&[{\cal D}_{\l},a_{\m}] = (h-\l s){\cal D}_{\l}c_{\m} -\{(h+\m
s)a_{\m}d_{\l} -
(h+\l s)c_{\m}b_{\l} + f(\m,\l)c_{\m}d_{\l}\}c_{\l}\cr
&&\qquad \qquad +\{(h+\l s)a_{\l}c_{\m} - (h+ \m s)c_{\l}a_{\m} +
f(\l,\m)c_{\l}c_{\m}\}\{(h+\l s)c_{\l} - d_{\l}\},\cr
&&[{\cal D}_{\l},b_{\m}] = (h+\l s)a_{\m}{\cal D}_{\l} + (h-\l
s){\cal D}_{\l}d_{\m}
+s(\m-\l)(h+\l s)c_{\l}a_{\m}c_{\l} \cr
&&\qquad \qquad +\{(h+\l s)a_{\l}d_{\m} - (h+\m s)c_{\l}b_{\m}
+ f(\l,\m)c_{\l}d_{\m}\}\{(h+\l s)c_{\l} - d_{\l}\} \cr
&&\qquad \qquad -a_{\l}\{(h-\l s)a_{\m}d_{\l} - (h-\m s)b_{\m}c_{\l}
- f(\l,\m)a_{\m}c_{\l}\} \cr
&&\qquad \qquad -(h+\l s)\{(h+\m s)a_{\m}d_{\l} - (h+\l s)c_{\l}b_{\m}
+ f(\m,\l)c_{\m}d_{\l}\}c_{\l}, \cr
&&[{\cal D}_{\l},c_{\m}] = 0 , \cr
&&[{\cal D}_{\l},d_{\m}] = (h+\l s){\cal D}_{\l}c_{\m} + s(\m-\l)(h+\l
s)a_{\l}c_{\m}c_{\l} +2hs(\l-\m)c_{\m}a_{\l}c_{\l} \cr
&&\qquad \qquad - a_{\l}\{(h-\l s)c_{\m}d_{\l} - (h-\m s)d_{\m}c_{\l}
-f(\l,\m)c_{\m}c_{\l}\} \cr
&&\qquad \qquad -\{(h-\m s)a_{\l}d_{\m} - (h-\l s)c_{\m}b_{\l}
+ (h-\m s)^2c_{\l}d_{\m}\}c_{\l}.
\eea
Clearly the determinant is not central unless the colour parameters
vanish.
These relations also exhibit an interesting feature that the determinants
associated with different colours are of  non-commuting nature, i.e.
$[{\cal D}_{\l},{\cal D}_{\m}] \ne 0$, unlike the coloured $q$-deformed
$GL(2)$ where they behave as commuting operators \cite{basu}.

If ${\cal D}_{\l}$ is non-singular for every $\l$, then the above
structure can be extended to a Hopf algebra by suitably defining
left and right inverses as
\bea
{\gamma }(T_{\l})&=&{{\cal D}_{\l}}^{-1}
\left(\begin{array}{cc}
 d_{\l}-(h-\l s)c_{\l} & -b_{\l}-(h-\l s)(d_{\l}-a_{\l})+(h-\l s)^2c_{\l}
\\
- c_{\l} & a_{\l} +(h-\l s)c_{\l} \end{array}\right)  \cr
&~& \cr
&=&\left(\begin{array}{cc}
 d_{\l}-(h+\l s)c_{\l} & -b_{\l}-(h+\l s)(d_{\l}-a_{\l})+(h+\l s)^2c_{\l}
\\
- c_{\l} & a_{\l} +(h+\l s)c_{\l} \end{array}\right) {{\cal
D}_{\l}}^{-1}\label{bdd}
\eea

The contents of this section generalize the results of \cite{aga}.
The entire Hopf algebra structure of $GL_{z,z^{\prime}}(2)$
can be easily recovered in the monochromatic limit $\l = \m$
and letting $(z, z^{\prime})$ in (\ref{ar}) go to $(-z, -z^{\prime})$.

\sect{Conclusions}

In this work, we have constructed a new
two-parametric deformation
of the Lie algebra $gl(2)$ . It is shown that this quantum algebra is
triangular (Jordanian). Starting
from the colour representation of its universal $R$-matrix, we have
obtained a new, non-standard
solution of the coloured YBE for $GL(2)$. Using this coloured solution we
have also
explored the quantum group structure associated with it.
We have thus obtained a coloured Jordanian quantum group
$GL_{h,s}^{(\l,\m)}(2)$ (infinite dimensional) incorporating the one
parameter $GL_z(2)$ and the two
parameter $GL_{z,z^{\prime}}(2)$ Jordanian quantum groups as
its finite-dimensional subalgebras.
It would be interesting to find out if the Jordanian $R$-matrix and its
super extension given in \cite{dp} can be regarded as different
manifestations of a  more general coloured solution.

{\bf Acknowledgements}

I would like to thank Prof. Edward Letzter for warm hospitality
at Texas A \& M University.

\bigskip



\end{document}